\documentclass{MBE}
\usepackage{amsmath}
 \usepackage{paralist}
  \usepackage{graphics}
   \usepackage{epsfig}
    \usepackage[colorlinks=true]{hyperref}
        \usepackage{ulem}

 \hypersetup{urlcolor=blue, citecolor=red}
\def\DOI{10.3934/mbe.2011.8.503}
\def\currentvolume{8}
 \def\currentissue{2}
  \def\currentyear{2011}
   \def\currentmonth{April}
    \def\ppages{503--513}
\newtheorem{theorem}{Theorem}[section]

\newtheorem{lemma}[theorem]{Lemma}
\newtheorem{proposition}[theorem]{Proposition}

\theoremstyle{definition}

\newtheorem{remark}{Remark}

\newcommand{\wA}{\widetilde{A}}
\newcommand{\ud}{\mathrm{d}}

\newcommand{\I}{\mathcal{I}}
\newcommand{\X}{\mathcal{X}}
\newcommand{\T}{\mathcal{T}}

\newcommand{\Y}{\mathcal{Y}}
\newcommand{\RR}{\mathbb{R}}

\title[Populations with diffusion and dynamic boundary
conditions]{Physiologically structured populations with diffusion and dynamic
boundary conditions}

\author[J\'{o}zsef Z. Farkas and Peter Hinow]{}
\email{jzf@maths.stir.ac.uk}
\email{hinow@uwm.edu}

\begin{document}
\maketitle

\centerline{\scshape J\'{o}zsef Z. Farkas }
\medskip
{\footnotesize
\centerline{Department of Computing Science and
Mathematics} \centerline{University of Stirling}
\centerline{Stirling, FK9 4LA, United Kingdom}
}
\medskip

\centerline{\scshape Peter Hinow}
\medskip
 {\footnotesize
\centerline{Department of
Mathematical Sciences}
\centerline{University of Wisconsin -- Milwaukee}
\centerline{P.O.~Box 413, Milwaukee, WI 53201, USA}}
%\bigskip
%
%\centerline{(Communicated by Nir Sochen)}

\begin{abstract}
We consider a linear size-structured population model with diffusion in the
size-space. Individuals are recruited into the population at arbitrary sizes.
We equip the model with generalized Wentzell-Robin (or dynamic) boundary
conditions. This approach allows the modelling of populations in which
individuals may have distinguished physiological states. We establish
existence and positivity of solutions by showing that solutions are governed by
a positive quasicontractive semigroup of linear o\-pe\-ra\-tors on the
biologically relevant state space. These results are obtained by establishing
dissipativity of a suitably perturbed semigroup generator. We also show that
solutions of the model exhibit balanced exponential growth, that is, our model
admits a finite-dimensional global attractor. In case of strictly positive
fertility we are able to establish that solutions in fact exhibit asynchronous
exponential growth.
\end{abstract}

\section{Introduction}

A significant amount of interest has been devoted to the analysis of
mathematical models arising in structured population dynamics (see
e.g.~\cite{MD,WEB} for references). Such models often assume spatial homogeneity
of the population in a given habitat and only focus on the  dynamics of the
population arising from differences between individuals with respect to some
physiological structure. In this context, reproduction, death and growth
characterize individual behavior which may be affected by competition, for
example for available resources.

In a recent paper, Hadeler \cite{H}  introduced size-structured population
models with diffusion in the size-space. The biological motivation is that
diffusion allows for ``stochastic noise'' to be incorporated in the model
equations in a deterministic fashion. The main question addressed in \cite{H}
is what type of boundary conditions are necessary for a biologically plausible
and mathematically sound model. In this context some special cases of a general
Robin boundary condition were considered. Diffusion terms have been introduced
into  structured population models by Waldst\"atter \textit{et
al.}~\cite{WHG}, Milner and Patton \cite{MP03} and Langlais and Milner
\cite{LM03} in the context of host-parasite models, where the parasite load is
the continuous structure variable. For example, when the
structuring variable represents a parasite load as in \cite{WHG}, at least one
special compartment of individuals arises, namely the class of the uninfected
ones. In a model where the structuring variable is continuous this state
corresponds to a set of measure zero, hence it does not carry mass. The true
meaning and advantage of employing Wentzell boundary conditions is that it
allows this special state to carry mass.

In this paper we introduce the following linear size-structured population
model
\begin{align}
&u_t(s,t)+
\left(\gamma(s)u(s,t)\right)_s\nonumber\\
=&\left(d(s)u_s(s,t)\right)_s-\mu(s)u(s,t)
+\int_0^m\beta(s,y)u(y,t)\,\ud y,\quad
s\in(0,m),\label{equation1} \\
& \left[\left(d(s)u_{s}(s,t)\right)_s\right]_{s=0}-b_0u_s(0,t)+c_0u(0,t)=0,
\label{equation2} \\
&\!\!\!
\left[\left(d(s)u_{s}(s,t)\right)_s\right]_{s=m}+b_mu_s(m,t)+c_mu(m,t)=0,
\label{equation3}
\end{align}
with a suitable initial condition.
The function $u=u(s,t)$ denotes the density of individuals of size, or
other developmental stage, $s$ at time $t$. Note that, we use $0$ as the minimal
value of the variable of $s$ only for mathematical convenience. It may well be
replaced by an arbitrary value $s_{min}$. The non-local integral term in
\eqref{equation1} represents the recruitment of individuals into the population.
Individuals may have different sizes at birth and $\beta(s,y)$ denotes the rate
at which individuals of size $y$ ``produce'' individuals of size $s$. Further
biologically relevant assumptions may be made on the fertility, such as
$\beta(s,y)=0$ for $s\ge y$, i.e.~parents cannot have larger offspring. We also
note that from this general model a single state at birth model may be deduced
by formally replacing the fertility function $\beta$ with an appropriate delta
function, see e.g.~\cite{MD}; Chapter I Section 4. $\mu$ denotes the
size-specific mortality rate while $\gamma$ denotes the growth rate. $d$ stands
for the size-specific diffusion coefficient, which we assume to be strictly
positive. $b_0$ and $b_m$ are positive numbers, while $c_0$ and $c_m$ are
non-negative. We will discuss special values for these constants in Equation
\eqref{condonbandc} below.

Equation \eqref{equation1} describes the evolution of a ``proper''
size-structured population, in contrast to the one where it is assumed that
all newborns enter the population at a minimal size. In that case, assuming
that the growth rate is positive, i.e.~individuals do not shrink, the linear
model can be rewritten as an age-structured model, see e.g.~\cite{MD}.
We also refer the interested reader to \cite{CS2,FGH,FHin} where size-structured
models with distributed recruitment processes were investigated. We also note
that although we refer to the structuring variable $s$ as size, it could well
represent any other physiological characteristic of individuals such as
accumulated energy or biomass, volume etc. The boundary conditions
\eqref{equation2}-\eqref{equation3} are the so called generalized
Wentzell-Robin or dynamic boundary conditions. These ``unusual'' boundary
conditions were investigated recently for models describing physical processes
such as diffusion and wave propagation, see e.g.~\cite{FGGR,FGGR2,GG}. Briefly,
they are used to model processes where particles reaching the boundary of a
domain can be either reflected from the boundary or they can be absorbed. Hence
the boundary points can carry mass. Our goal here is to introduce this
rather general type of boundary condition in the context of models describing
the evolution of biological populations, with particular focus on positivity and
asymptotic behavior of solutions. Potential applications include cell
populations with resting states at $s=0$ and $s=m$, or models for populations
structured by an infection level or parasite load \cite{WHG}.

The first works introducing boundary conditions that involve second order
derivatives for parabolic or elliptic differential operators go back to the
1950s, see  the papers by Feller \cite{F1,F2} and Wentzell \cite{V1,V2}.
These first studies were purely motivated from the mathematical point of view.
The original question was to identify the set of all possible boundary
conditions that give rise for a parabolic differential operator to generate a
contraction semigroup on an appropriate state space. The abstract mathematical
analysis gave a clue where physical intuition failed previously.

The boundary conditions \eqref{equation2}-\eqref{equation3} are in general form,
and we shall now specify the constants $b_0,b_m$ and $c_0,c_m$ to give a
biological explanation for the boundary conditions. Integration of Equation
\eqref{equation1} from $0$ to $m$ yields for
$U(t)=\int_0^m u(s,t)\,\ud s$
\begin{equation*}
\begin{aligned}
\frac{d}{dt}\,U(t) & = \gamma(0)u(0,t)-\gamma(m)u(m,t)
+d(m)u_s(m,t)-d(0)u_s(0,t) \\
&\phantom{=} +\int_0^m\int_0^m\beta(s,y)u(y,t)\,\ud y\,\ud
s-\int_0^m\mu(s)u(s,t)\,\ud s \\
&=:B(t)-D(t)+\gamma(0)u(0,t)-\gamma(m)u(m,t) +d(m)u_s(m,t)-d(0)u_s(0,t),
\end{aligned}
\end{equation*}
where $B$ and $D$ denote the combined birth and death processes, respectively.
We note that formally, by replacing the diffusion term by its counterpart from
Equation \eqref{equation1}, the boundary conditions
\eqref{equation2}-\eqref{equation3} can be cast in the {\it dynamic form}
\begin{equation}\label{dynamicbc}
\begin{aligned}
u_t(0,t) =& u(0,t)(-\gamma'(0)-\mu(0)-c_0)+u_s(0,t)(b_0-\gamma(0))
+\int_0^m\beta(0,y)u(y,t)\,\ud y, \\
u_t(m,t) =&
u(m,t)(-\gamma'(m)-\mu(m)-c_m)+u_s(m,t)\left(-b_m-\gamma(m)\right)
\\&+\int_0^m\beta(m,y)u(y,t)\,\ud y.
\end{aligned}
\end{equation}
These are the governing equations for individuals of minimum and maximum sizes,
respectively. It is natural to assume that in the absence of mortality and
recruitment, i.e.~when $B(\cdot)\equiv D(\cdot)\equiv 0$, the total population
size $U(t)+u(0,t)+u(m,t)$ remains constant at every time $t$.
Mathematically, this amounts to the
condition
\begin{equation*}
\begin{aligned}
0=&\frac{d}{dt}U(t)+u_t(0,t)+u_t(m,t)\\=&u(0,
t)\left(\gamma(0)-\gamma'(0)-c_0\right)+u(m,t)
\left(-\gamma(m)-c_m-\gamma'(m)\right)\\
&\phantom{=}+u_s(0,t)\left(b_0-d(0)-\gamma(0)\right)+u_s(m,
t)\left(d(m)-\gamma(m)-b_m\right).
\end{aligned}
\end{equation*}
In order to guarantee conservation of total population in the absence of birth
and death processes, we make the following assumptions
\begin{equation}\label{condonbandc}
c_0=\gamma(0)-\gamma'(0),\,c_m=-\gamma(m)-\gamma'(m),\,b_0=d(0)+\gamma(0),\,
b_m=d(m)-\gamma(m).
\end{equation}
We note that condition \eqref{condonbandc} together with the assumption that $b$
is positive and $c$ is non-negative impose a restriction on the growth rate $\gamma$.

The dynamic boundary conditions \eqref{dynamicbc} can now be written as
\begin{equation}\label{dynamicbc2}
\begin{aligned}
u_t(0,t) & = (-\mu(0)-\gamma(0))u(0,t)+d(0)u_s(0,t)+\int_0^m\beta(0,y)u(y,t)\,\ud y, \\
u_t(m,t) & =(\gamma(m)-\mu(m))u(m,t)-d(m)u_s(m,t)+\int_0^m\beta(m,y)u(y,t)\,\ud y.
\end{aligned}
\end{equation}
Hence the dynamics of individuals in the two special states $s=0$ and $s=m$ are
governed by equations \eqref{dynamicbc2}. The meaning of the governing equations
\eqref{dynamicbc2} is intuitively clear. For example, at $s=0$, individuals are
leaving this compartment due to growth and mortality and are recruited according
to the integral term, while the diffusion accounts for the flux through this
state (this also yields a loss term, since the outer normal derivative is
$-u_s$).

We impose the following assumptions on the model ingredients
\begin{equation*}
\mu\in C([0,m]),\quad \beta\in C([0,m]\times [0,m]),\quad \beta,\,\mu \geq
0,\quad \gamma, d\in C^1([0,m]),\quad d>0.
\end{equation*}
In this note first we establish existence and positivity of solutions of model
\eqref{equation1}-\eqref{equation3}. This existence proof follows similar
arguments developed in \cite{FGGR,FGGR2}. The significant difference is that
solutions to our model are not necessarily governed by a contraction
semigroup, hence as in \cite{FH2,FH3} we need to rescale the semigroup to obtain
the dissipativity estimate. As a result of the dissipativity
calculation,
we show that the resolvent operator of the semigroup generator is
positive. This was not established in \cite{FGGR2}.

In Section \ref{section:asybehavior} we investigate the asymptotic behavior of
solutions. First we establish that solutions of the model equations exhibit
balanced exponential growth, in general. This is an interesting phenomenon,
often observed for linear structured population models, see
e.g.~\cite{FH2,FHin}. In some sense, it is a stability result and characterizes
the global asymptotic behavior of solutions to the model. Then, assuming that
fertility is strictly positive, we are able to show that after a rescaling by
and exponential factor solutions actually tend to a fixed size-distribution.
This is shown via establishing irreducibility of the governing semigroup.

\section{Existence and positivity of solutions}\label{section:existence}
In this section we are going to establish the existence of a positive
quasicontractive semigroup of operators which governs the evolution of solutions
of \eqref{equation1}-\eqref{equation3}.
For basic definitions and results used throughout this paper we refer the reader
to \cite{Adams,AGG,CH,NAG}.
Let\\
 $$\mathcal{X}=\left(L^1(0,m)\oplus\mathbb{R}^2,||\cdot||_{\X}\right),$$

 \noindent where
for \mbox{$(x,x_0,x_m)\in\mathcal{X}$} the norm is given by
\begin{equation}\label{choice_of_norm}
||(x,x_0,x_m)||_{\mathcal{X}}=||x||_{L^1}+c_1|x_0|+c_2|x_m|,
\end{equation}
for some positive constants $c_1$ and $c_2$ that we will specify later.
Then $\mathcal{X}$ is a Banach lattice.
We identify  a function $u\in C[0,m]$ with its restriction triple
$(u|_{(0,m)},u(0),u(m))\in \mathcal{X}$. With this identification, the set
$C^2[0,m]$ is dense in $\mathcal{X}$ with respect to the
$||\,\cdot\,||_{\mathcal{X}}$-norm. Let
\begin{equation*}
\begin{aligned}
D(A) &= \left\{ u\in C^2[0,m] \: :\: \Psi u \in L^1(0,m),\quad
(d(s)u'(s))'\big|_{s=0}-b_0u'(0)+c_0u(0)=0,   \right. \\
     &\phantom{=} \:\,\left. (d(s)u'(s))'\big|_{s=m}+b_mu'(m)+c_mu(m) =0
\right\},
\end{aligned}
\end{equation*}
where
\begin{equation*}
\Psi u(s) = (d(s)u'(s))' -(\gamma(s)u(s))'
-\mu(s)u(s)+\int_0^m\beta(s,y)u(y)\,\ud y.
\end{equation*}
The operator $A$ with domain $D(A)$ is then
defined by
\begin{equation*}
Au = \begin{pmatrix}
\Psi u  \vspace{1mm} \\
(b_0-\gamma(0))\frac{d}{ds}u\,\big|_{s=0}+\int_0^m\beta(0,y)u(y)\,\ud
y-\rho_0u(0)  \vspace{1mm} \\
(-b_m-\gamma(m))\frac{d}{ds}u\,\big|_{s=m}+\int_0^m\beta(m,y)u(y)\,\ud
y-\rho_mu(m)\\
\end{pmatrix}
\end{equation*}
where we have set
\begin{equation*}
\rho_0 = \mu(0)+c_0+\gamma'(0),\quad \rho_m = \mu(m)+c_m+\gamma'(m),
\end{equation*}
for short. We use \cite[Theorem 3.15]{NAG} (see also \cite[Section A-II, Theorem
2.11]{AGG}) which characterizes generators of contractive semigroups via
dissipativity. We establish the existence of a quasicontractive semigroup for
the general boundary condition \eqref{equation2}-\eqref{equation3}. Recall that
a strongly continuous semigroup
$\left\{\mathcal{T}(t)\right\}_{t\ge 0}$ is called \textit{quasicontractive} if
$||\T(t)||\le e^{\omega t}$ for some $\omega\in\RR$ and it is  called
\textit{contractive} if the choice  $\omega\le 0$ is possible. A linear
operator $A$ with domain $D(A)$ is  \textit{dissipative}, if for all $x\in
D(A)$ and $\lambda>0$ one has $||(\lambda\I-A)x||\ge \lambda||x||$.

\begin{theorem}\label{existence_theorem} Assume that
\begin{equation}\label{theoremcond}
c_1=\frac{d(0)}{b_0-\gamma(0)}>0,\quad \text{and}\quad
c_2=\frac{d(m)}{\gamma(m)+b_m}>0.
\end{equation}
Then the closure of the operator $A$ is the infinitesimal
generator of a positive quasicontractive semigroup of bounded linear
operators on the state space $\X$ (where the weights in equation
\eqref{choice_of_norm} are chosen accordingly).
\end{theorem}
\noindent  \emph{Proof.}
We introduce the modified operator $\wA$  with $\beta=0$ in the definition of
$A$. For $\lambda>0$ and $h\in \X, u\in D(\wA)$ we consider the
equation
\begin{equation}\label{resolventeq}
u-\lambda\left(\wA-\omega\mathcal{I}\right)u=h,\quad
\text{on}\quad [0,m].
\end{equation}
That is
\begin{align}
h(s) & =u(s)-\lambda\left(\left(d(s)u'(s)\right)'-
\left(\gamma(s)u(s)\right)'-\mu(s)u(s)
-\omega u(s)\right),\quad s\in (0,m), \label{resolventeq2} \\
\lambda^{-1}h(0) & =
u'(0)(\gamma(0)-b_0)+u(0)\left(\lambda^{-1}+\gamma'(0)+\mu(0)+c_0+\omega
\right),
\label{resboundc3} \\
\lambda^{-1}h(m) & =
u'(m)(\gamma(m)+b_m)+u(m)\left(\lambda^{-1}+\gamma'(m)+\mu(m)+c_m+\omega
\right).\label{resboundc4}
\end{align}
Next we multiply Equation \eqref{resolventeq2} by $\chi_{u^+}(s)$, where
$\chi_{u^+}$ denotes the characteristic function of $u^+$,
and integrate from $0$ to $m$. The boundary $\partial [u>0]$
consists of two parts, $\Gamma_1 = \partial [u>0]\cap (0,m)$ and
$\Gamma_2=\overline{[u>0]}\cap \{0,m\}$. The term $d(s)u'(s)$ gives negative
contributions on $\Gamma_1$ since the outer normal derivative of $u$ is
negative.  The term $\gamma(s)u(s)$ gives no contributions on $\Gamma_1$ since
$u=0$ there. Hence we obtain
\begin{align*}
||u^+||_1\le &\lambda\left(\text{sgn}(u^+(m))d(m)u'(m)
-\text{sgn}(u^+(0))d(0)u'(0)\right)\\
& -\lambda\left(\text{sgn}(u^+(m))\gamma(m)u(m)-\text{sgn}
(u^+(0))\gamma(0)u(0)\right) \\
& -\lambda\int_0^m\chi_{u^+}(s)\mu(s)u(s)\,\ud s-\lambda \omega
||u^+||_1+\int_0^m\chi_{u^+}(s)h(s)\,\ud s.
\end{align*}
This, combined with equations \eqref{resboundc3}-\eqref{resboundc4}, yields
\begin{equation}\label{pluspart}
\begin{aligned}
& ||u^+||_1 +\text{sgn}(u^+(m))u(m)\left(\frac{d(m)}{\gamma(m)+b_m}
\left(1+\lambda\left(\gamma'(m)+\mu(m)+c_m+\omega\right)\right)
+\lambda\gamma(m)\right)\\
&\quad+\text{sgn}(u^+(0))u(0)\left(\frac{d(0)}{b_0-\gamma(0)}
\left(1+\lambda\left(\gamma'(0)+\mu(0)+c_0+\omega\right)\right)
-\lambda\gamma(0)\right) \\
& \le -\lambda\int_0^m\chi_{u^+}(s)\mu(s)u(s)\,\ud s-\lambda
\omega ||u^+||_1+\int_0^m\chi_{u^+}(s)h(s)\,\ud s \\
&\quad+\text{sgn}(u^+(m))h(m)\frac{d(m)}{\gamma(m)+b_m}+\text{sgn}
(u^+(0))h(0)\frac{d(0)}{b_0-\gamma(0)}.
\end{aligned}
\end{equation}
Similarly, multiplying Equation \eqref{resolventeq2} by $-\chi_{u^-}(s)$ and
integrating from $0$ to $m$, we obtain
\begin{equation}\label{negpart}
\begin{aligned}
& ||u^-||_1 -\text{sgn}(u^-(m))u(m)\left(\frac{d(m)}{\gamma(m)+b_m}
\left(1+\lambda\left(\gamma'(m)+\mu(m)+c_m+\omega\right)\right)
+\lambda\gamma(m)\right) \\
& \quad-\text{sgn}(u^-(0))u(0)\left(\frac{d(0)}{b_0-\gamma(0)}
\left(1+\lambda\left(\gamma'(0)+\mu(0)+c_0+\omega\right)\right)
-\lambda\gamma(0)\right) \\
&\le\lambda\int_0^m\chi_{u^-}(s)\mu(s)u(s)\,\ud s-\lambda\omega
||u^-||_1-\int_0^m\chi_{u^-}(s)h(s)\,\ud s \\
& \quad-\text{sgn}(u^-(m))h(m)\frac{d(m)}{\gamma(m)+b_m}-\text{sgn}
(u^-(0))h(0)\frac{d(0)}{b_0-\gamma(0)}.
\end{aligned}
\end{equation}
Adding \eqref{pluspart} and \eqref{negpart} yields
\begin{equation}\label{sumpart}
\begin{aligned}
&||u||_1+|u(0)|\left(\frac{d(0)}{b_0-\gamma(0)}\left(1+\lambda
\left(\gamma'(0)+\mu(0)+c_0+\omega\right)\right)-\lambda\gamma(0)\right)\\
&\quad+|u(m)|\left(\frac{d(m)}{\gamma(m)+b_m}\left(1+\lambda
\left(\gamma'(m)+\mu(m)+c_m+\omega\right)\right)+\lambda\gamma(m)\right)\\
\le&\int_0^m\left(\chi_{u^+}(s)-\chi_{u^-}(s)\right)h(s)\,\ud
s+h(0)\frac{d(0)}{b_0-\gamma(0)}\left(\text{sgn}(u^+(0))-\text{sgn}
(u^-(0))\right)\\&-\lambda \omega ||u||_1
+h(m)\frac{d(m)}{\gamma(m)+b_m}\left(\text{sgn}(u^+(m))-\text{sgn}
(u^-(m))\right)\\&-\lambda\int_0^m\left(\chi_{u^+}(s)-\chi_{u^-}(s)\right)
\mu(s)u(s)\,\ud s.
\end{aligned}
\end{equation}
Assuming that condition \eqref{theoremcond} holds true the left hand side
of inequality \eqref{sumpart} can be estimated below, for $\omega\in\mathbb{R}$
large enough, by
\begin{equation*}
||u||_1+|u(0)|\frac{d(0)}{b_0-\gamma(0)}+|u(m)|\frac{d(m)}{\gamma(m)+b_m}.
\end{equation*}
Similarly, the right hand side of inequality \eqref{sumpart} can be estimated
above by
\begin{equation*}
||h||_1+|h(0)|\frac{d(0)}{b_0-\gamma(0)}+|h(m)|
\frac{d(m)}{\gamma(m)+b_m}.
\end{equation*}
Hence, if condition \eqref{theoremcond} is satisfied  we have the dissipativity
estimate
\begin{equation*}
||u||_{\X}\le||h||_{\X}=\left|\left|u-\lambda\left(\wA-\omega\I\right)u
\right|\right|_{\X}
\end{equation*}
for the operator $\wA-\omega\I$.

For the range condition we need to show that whenever $h$ is in a dense
subset of $\mathcal{X}$ then the solution $u$ of Equation \eqref{resolventeq}
belongs to the domain of $A$. Since we assumed that $d>0$, i.e.~we have true
Wentzell boundary conditions at both endpoints of the domain, the required
regularity of the solution $u\in C^2[0,m]$ and hence the range condition follows
from \cite[Theorem 6.31]{GT}. Thus the closure of $\wA-\omega\I$ is a generator
of a contractive semigroup by \cite[Theorem 3.15]{NAG}.

Next we observe that for $h\in\mathcal{X}^+$ every term on the right hand
side of Equation \eqref{negpart} is non-positive, while every term on the left
hand side of Equation \eqref{negpart} is non-negative. The inequality can only
hold for $h\in\mathcal{X}^+$ if $||u^-||_1=u(0)=u(m)=0$.
This proves that the resolvent operator
$R(\lambda,\wA-\omega\I)=(\lambda\mathcal{I}-(\wA-\omega\I))^{-1}$ is positive
for $\lambda$
large enough, hence the closure of $\wA-\omega \mathcal{I}$ generates a
positive semigroup. A simple perturbation result yields that
the closure of $\wA$ is a generator of a positive
quasicontractive semigroup.

Finally we note that, the operator $A-\wA$ is positive, bounded and linear,
hence it generates a positive semigroup
$\left\{\mathcal{S}(t)\right\}_{t\ge 0}$ which
satisfies
\begin{equation*}
||\mathcal{S}(t)||\le e^{tB}, \quad \text{where}\quad B=||\beta||_\infty.
\end{equation*}
The proof of the Theorem is now completed on the grounds of the Trotter product
formula, see e.g.~\cite[Corollary 5.8 Ch. 3]{NAG}. \hfill $\Box$
\begin{remark} Conditions \eqref{theoremcond} are clearly satisfied if only the
denominators are positive. However, the notation in \eqref{theoremcond} gives
immediately the norm on the state space.
\end{remark}
\begin{remark} We note that the biologically natural assumptions in
\eqref{condonbandc} that guarantee conservation of the total population
in the absence of mortality and recruitment imply that $c_1=c_2=1$. Hence the
mathematical calculations coincide with the biological intuition.
It can be shown that in the absence of mortality and recruitment the
operator $A$ generates a contraction semigroup on the state space
$\left(\X,||\cdot||_{\X}\right)$, in fact
$||\T(t)||_{\mathcal{X}}=1$ for all $t\ge0$.
\end{remark}

\section{Asymptotic behavior}\label{section:asybehavior}
In this section we investigate in the framework of semigroup theory the
asymptotic behavior of solutions of model
\eqref{equation1}-\eqref{equation3}. Since our model is a linear one, we may
expect that solutions either grow or decay exponentially, unless they are
at a steady state. For simple structured population models in fact it is often
observed (see e.g.~\cite{DVBW,FH2,FHin,GW}), that solutions grow exponentially
and tend to a finite-dimensional (resp.~one-dimensional) global
attractor. This phenomenon is called \textit{balanced}
(resp.~\textit{asynchronous}) \textit{exponential growth}.
The rate of exponential growth is called the \textit{Malthusian parameter} or
\textit{intrinsic growth rate}. In other words, asynchronous exponential growth
means that solutions tend to a fixed size distribution (often called the stable
size profile) after a suitable rescaling of the semigroup.

In this section we show that solutions of our model exhibit the same asymptotic
behavior. This question can be addressed effectively in the framework of
semigroup theory (see \cite{AGG,CH,NAG}). Briefly, to establish balanced
exponential growth, one needs to show that the growth bound of the semigroup is
governed by a leading eigenvalue of finite (algebraic) multiplicity of its
generator, and there exists a spectral gap, i.e.~the leading eigenvalue is
isolated in the spectrum. Our first result will assure this latter condition.
Moreover, if it is also possible to establish that the semigroup is irreducible
then one has that the algebraic multiplicity of the spectral bound equals one
(with a corresponding positive eigenvector) and the semigroup exhibits
asynchronous exponential growth.
\begin{lemma}
The spectrum of $A$ can contain only isolated eigenvalues of finite
algebraic multiplicity.
\end{lemma}
\noindent \emph{Proof.} We show that the resolvent operator $R(\lambda,A)$ is
compact. Since $\widetilde{A}$ is a bounded perturbation of $A$ it is enough to
show that $R(\lambda,\widetilde{A})$ is compact. This follows however, from the
regularity of the solution of the resolvent Equation \eqref{resolventeq}
and noting that $C^2[0,m]\subset W^{1,1}(0,m)\oplus\mathbb{R}^2$ which is
compactly embedded in $L^1(0,m)\oplus\mathbb{R}^2$
by the Rellich-Kondrachov Theorem \cite[Theorem 6.3, Part I]{Adams}. The claim
follows on the ground of \cite[Proposition II.4.25]{NAG}.
\hfill $\Box$

Next we recall some necessary basic notions from linear semigroup theory, see
e.g.~\cite{AGG,CH,NAG}.
A strongly continuous semigroup $\left\{\mathcal{S}(t)\right\}_{t\geq 0}$ on a
Banach space $\Y$ with generator
${\mathcal O}$ and \textit{spectral bound}
\begin{equation*}
s\left({\mathcal O}\right)=\sup
\left\{\,\text{Re}(\lambda)\::\:\lambda\in\sigma\left({\mathcal
O}\right)\,\right\}
\end{equation*}
is said to exhibit \textit{balanced exponential growth} if there exists a
projection $\Pi$ on $\Y$ such that
\begin{equation*}
\lim_{t\to \infty} \|e^{-s\left({\mathcal O}\right)t}\,\mathcal{S}(t)-\Pi\|=0.
\end{equation*}
If the projection $P$ is of rank one then the semigroup has {\it asynchronous
exponential growth}.
Moreover, the \textit{growth bound} $\omega_0$ is the infimum of all real
numbers $\omega$ such that there exists a constant $M\ge 1$ with
$||\mathcal{S}(t)||\le Me^{\omega t}$. We also recall (see e.g.~\cite[C-III
Definition 3.1]{AGG}) that a positive semigroup
$\left\{\mathcal{S}(t)\right\}_{t\ge 0}$ on a Banach lattice $\mathcal{Y}$ is
called \textit{irreducible} if there is no $\mathcal{S}(t)$ invariant
closed ideal of $\mathcal{Y}$ except the trivial ones, $\{0\}$ and
$\mathcal{Y}$.
\begin{theorem}
Model \eqref{equation1}-\eqref{equation3} admits a finite dimensional global
attractor.
\end{theorem}
\noindent \emph{Proof.} We have shown in Theorem \ref{existence_theorem} that
solutions of our model are governed by a positive semigroup. Derndinger's
Theorem (see e.g.~\cite[Theorem VI.1.15]{NAG})
implies that the spectral bound of the generator equals to the growth bound of
the semigroup, i.e.~$s(A)=\omega_0$. Lemma 3.1 implies that the spectral bound
$s(A)$ is an eigenvalue of finite algebraic (hence geometric) multiplicity
unless the spectrum is empty. If the spectrum of $A$ is empty we have by
definition $\omega_0=-\infty$ and every solution tends to zero. Otherwise we
have for the semigroup $\left\{\mathcal{T}(t)\right\}_{t\ge 0}$ generated by the
closure of $A$
\begin{equation}
\lim_{t\to\infty}\left|\left|e^{-s(A)t}\mathcal{T}(t)-\Pi\right|\right|=0,
\end{equation}
where $\Pi$ is the projection onto the finite dimensional eigenspace
corresponding
to the eigenvalue $s(A)$.
\hfill $\Box$

Recall that a subspace $I\subset\mathcal{Z}$ of a Banach
lattice $\mathcal{Z}$ is an \textit{ideal} iff $f\in I$ and $|g|\le |f|$
implies that $g\in I$. The following useful result is from \cite[C-III,
Proposition 3.3]{AGG}.
\begin{proposition}\label{idealchar}
Let $B$ be the generator of a positive semigroup $\mathcal{U}(t)$ on the
Banach lattice $\mathcal{Z}$ and $K$  a bounded positive operator.
Let $\mathcal{V}(t)$ be the  semigroup generated by $B+K$. For a closed ideal
$I\subset\mathcal{Z}$ the following assertions are equivalent:
\\  (i) $I$ is $\mathcal{V}$-invariant,
\\  (ii) $I$ is invariant under both $\mathcal{U}$ and $K$.
\end{proposition}
We introduce the recruitment operator $K=A-\widetilde{A}$.

\begin{theorem}\label{aeg}
Assume that \eqref{theoremcond} holds true and $\beta>0$. Then the semigroup
generated by the closure of $A$ exhibits
asynchronous exponential growth.
\end{theorem}
\noindent \emph{Proof.} Our goal is to apply Proposition \ref{idealchar} for the
operator $\widetilde{A}$,  whose closure is the generator  of a positive
semigroup as shown in Theorem \ref{existence_theorem} and to $K$,
which is clearly positive and bounded.
Every closed ideal $I$ of $\mathcal{X}$ can be written as $I_1\oplus I_2\oplus
I_3$, where $I_1$ is a closed ideal in
the Banach lattice $L^1(0,m)$ and $I_2,I_3$ are closed ideals in $\mathbb{R}$.
Note that, $\mathbb{R}$ admits only two ideals, i.e.~$\{0\}$ or $\mathbb{R}$
itself.
Next we observe that non-trivial closed ideals in $L^1(0,m)$ can be
characterized
via closed subsets $G$ of positive measure of $(0,m)$. That is,
the subspace $J$ is a closed ideal of $L^1(0,m)$ if it contains the functions
$f\in L^1(0,m)$ vanishing on $G$.
Next we show that no non-trivial closed ideal $I=I_1\oplus I_2\oplus I_3$ is
invariant under $K$ or under the semigroup generated by $\widetilde{A}$.
If $I_1\ne \{0\}$ then the condition $\beta>0$ guarantees that $K
u(s)=\int_0^m\beta(s,y)u(y)\,\ud y>0$, for every $s\in (0,m)$ for any $u\in
I_1$, i.e.~the image $Ku$ does not vanish anywhere, hence by the previous
characterization we must have $I_1=L^1(0,m)$. Moreover, in this case $\beta>0$
implies that we must have $I_2=I_3=\mathbb{R}$, since
$\int_0^m\beta(0,y)u(y)\,\ud y>0$ and $\int_0^m\beta(m,y)u(y)\,\ud y>0$ for any
$u\not\equiv 0$. On the other hand, if $I_1=\{0\}$ then we have
$D(\widetilde{A})\cap I=\{0\}$, hence the restriction of $\widetilde{A}$ to
$\mathbb{R}\oplus\mathbb{R}$ (or even to $\mathbb{R}$)  does not generate a
semigroup. This means that $I$ cannot be invariant under the semigroup generated
by $\widetilde{A}$. That is we have by Proposition \ref{idealchar} that the
semigroup generated by $A$ has no non-trivial closed invariant ideal. Therefore
it is irreducible, and solutions exhibit asynchronous exponential growth, see
e.g.~\cite{CH}.
\hfill $\Box$

The previous result completely characterizes the asymptotic behavior of
solutions to the population model. That is, solutions behave asymptotically as
$e^{rt}u_*(s)$ independently of the initial condition, where $r$ is the so
called Malthusian parameter, and $u_*$ is often referred to
as the final size distribution.

\section{Concluding remarks}

In this note we introduced a linear structured population model with diffusion
in the size space. Introduction of a diffusion is natural in the biological
context \cite{H}, since unlike in age-structured models, individuals that have
the same size initially, may disperse as time progresses.  In other words,
diffusion amounts to adding noise in a deterministic fashion. We equipped our
model with generalized Wentzell-Robin boundary conditions. We showed that the
model is governed by a positive quasicontractive semigroup on the biologically
relevant state space. Furthermore we have characterized the asymptotic behavior
of solutions via balanced exponential growth of the governing semigroup. We also
established that solutions exhibit asynchronous exponential growth if the
function $\beta$ is strictly positive. An important biological
consequence of asynchronous exponential growth is population stabilization in
the sense that the proportion of the population in any subset of the structure
space converges to a limiting value as time evolves, independently of the
initial state of the population. The question of irreducibility of the
semigroup generated by the Laplace operator with mixed Robin boundary conditions
on a $L^p$-space (for $1<p<\infty$) was addressed in the recent work by
Haller-Dintelmann \textit{et al.}~\cite{HDHR}. It is expected that a
similar result holds if generalized Wentzell-Robin boundary conditions are
imposed. This is a topic of ongoing research.

In our model we have taken the view that individuals may be recruited into the
population at different sizes. This appears to be the natural choice in the
context of general physiologically structured population models, as opposed to
age-structured models, where every individual is born at the same age zero. It
is interesting to investigate whether a mathematically sound ``limiting
relationship'' exists between models with infinite states at birth and one state
at birth. This will be addressed in future work.

The power of generalized Wentzell boundary conditions in the context of
population models is to allow the boundary states to carry mass. This is
especially interesting in the $L^p$ context as the boundary is a set of
measure zero and therefore seems to play no role in an integral term.
Interestingly, sinks on the boundary can cause ill-posedness in space dimensions
$\ge2$ as Vazquez and Vitillaro \cite{VV} have shown.

In the future we will extend our model to incorporate interaction variables,
to allow competition. Then model \eqref{equation1}-\eqref{equation3} becomes a
nonlinear one, and the mathematical analysis will become more  difficult. To our
knowledge, positivity results are rather rare in the literature for nonlinear
models. In \cite{FGGR} it was shown that the nonlinear semigroup generator
satisfies the positive minimum principle, hence the semigroup is positive. This
however, does not apply to population models, as the positive cone of the
natural state space $L^1$ has empty interior, hence the positive minimum
principle does not apply. It will be also interesting to consider more general
models with a finite number of structuring variables, such as age-size
structured models. Then the domain will be a cube  $[0,1]^n$ and the
prescription of appropriate boundary conditions will be much more involved.

\section*{Acknowledgments}
J.~Z.~Farkas was supported by a personal research grant from the Carnegie Trust
for the Universities of Scotland. Part of this work was done while J.~Z.~Farkas
visited the University of Wisconsin - Milwaukee. Financial support from the
Department of Mathematical Sciences is greatly appreciated. We thank two
anonymous referees for valuable comments and suggestions.

\quad Received March 4, 2010; Accepted December 3, 2010.

\end{document}